\documentclass[12pt]{article}

\usepackage{amsmath,amssymb,amsfonts,theorem,makeidx,latexsym,epsfig}

\usepackage{subfigure}

\newtheorem{defn}{Definition}[section]

\newtheorem{lemma}[defn]{Lemma}

{\theorembodyfont{\rmfamily}

\newtheorem{ex}[defn]{Example}}

\newtheorem{thm}[defn]{Theorem}

\newtheorem{prop}[defn]{Proposition}

\newtheorem{cor}[defn]{Corollary}

\newtheorem{rem}[defn]{Remark}

\numberwithin{equation}{section}

\newcommand{\ltr}{ L^2(\mathbb R) }

\newcommand{\mn}{\mathbb N}

\newcommand{\mr}{\mathbb R}

\newcommand{\mz}{\mathbb Z}

\newcommand{\mts}{ \{E_{mb}T_{na}g \}_{m,n \in \mz}}

\def\bp{{\noindent\bf Proof. \ }}

\def\ep{\hfill$\square$\par\bigskip}

\def\bqs{\begin{equation}}

\def\eqs{\tag*{$\square$}\end{equation}\par\bigskip}

\def\al{\alpha}

\def\supp{\text{supp}}

\def\vn{\vspace{.1in}\noindent}

\def\bop{\begin{op}\rm}

\def\eop{\end{op}}

\def\bee{\begin{eqnarray}}

\def\ene{\end{eqnarray}}

\def\bes{\begin{eqnarray*}}

\def\ens{\end{eqnarray*}}

\def\bei{\begin{itemize}}

\def\eni{\end{itemize}}

\def\bt{\begin{thm}}

\def\et{\end{thm}}

\def\bc{\begin{cor}}

\def\ec{\end{cor}}

\def\bpr{\begin{prop}}

\def\epr{\end{prop}}

\def\bl{\begin{lemma}}

\def\el{\end{lemma}}

\def\bd{\begin{defn}}

\def\ed{\end{defn}}

\def\bex{\begin{ex}}

\def\enx{\end{ex}}

\def\bfi{\begin{fig}}

\def\efi{\end{fig}}

{\nopagebreak\par\noindent\hbox to
\hsize{\hrulefill}\par\noindent}

\title{On Gabor frames generated by sign-changing windows and B-splines}

\date{\today}

\author{Ole Christensen\thanks{
Department of Applied Mathematics and Computer Science, Technical
University of Denmark, Building 303, 2800 Lyngby, Denmark.
E-mail: ochr@dtu.dk}, Hong Oh Kim\thanks{
Division of General Studies, UNIST, UNIST-gil 50, Ulsan 689-798, Republic of Korea.
E-mail: hkim2031@unist.ac.kr}, Rae Young Kim\thanks{
Department of Mathematics, Yeungnam University, 280 Daehak-Ro, Gyeongsan,
Gyeongbuk, 712-749, Republic of Korea.
E-mail: rykim@ynu.ac.kr
}}

\begin{document}

\maketitle

\begin{abstract} For a class of compactly supported windows we characterize the
frame property for a Gabor system  $\mts,$ for translation parameters $a$ belonging to  a certain range depending on the support size. We show that  the obstructions to the frame property are located on a countable number of ``curves."  For functions that are positive on the interior of the support these obstructions do not appear, and the considered region in the $(a,b)$ plane is fully contained in the frame set. In particular this confirms a recent conjecture about B-splines by Gr\"ochenig in that particular region.  We prove that the full conjecture is true if it can be proved in a certain ``hyperbolic strip."


\end{abstract}

\section{Introduction}

Only for quite special functions $g\in\ltr$ we know a characterization of
the {\it Gabor frame set,}
${\cal F}(g):= \{(a,b)\in \mr_+^2 \, \big| \, \mts \, \mbox{is a frame}\}; $ these functions include
the Gaussian \cite{Ly, Se2}, the hyperbolic secant \cite{JS2}, the one-sided/two-sided exponentials \cite{Jan9}, and totally positive functions \cite{GrSt}. Common for all these functions is that they are nonnegative.

Much less is known about more general functions, e.g., functions that change sign.
In this paper we consider a class of continuous compactly supported windows $g$ with $\supp \, g= [ - \alpha, \alpha]$ for some $\alpha>0$
and characterize the frame property of $\mts$ in the region
$ \alpha \le a <2\alpha, b< 1/a.$
For technical reasons (and in order to avoid pathological examples of no practical interest)
we assume that the function $g$ only has a finite number of zeros in $]-\alpha, \alpha[.$
The general result, to be stated in Theorem \ref{2-2-1}, shows that the zeros in the interior of the support
lead to certain obstacles for the frame property that cannot be predicted from the known results for nonnegative functions. For each translation parameter a countable number of obstructions
can appear, i.e., one can think about the obstructions as located on a countable number of
curves in the $(a,b)$--plane.
The general result also implies the existence of a compactly supported dual window
if the frame property is satisfied, with an interesting interpretation
in terms of the redundancy of the frame: in fact, if $\frac{M-1}{M} \le ab < \frac{M}{M+1}$ for some
$M=2,3,\dots,$ i.e.,
if the redundancy $(ab)^{-1}$ is
at least  $\frac{M+1}{M}= 1+ 1/M,$
the existence of a dual window supported on $[-2\alpha M, 2\alpha M]$ is guaranteed.

In the special case of a function $g$ that is positive on $]-\alpha, \alpha[$
the  general result implies that
$\mts$ is a frame for all parameters $a,b$ in the considered region
$\alpha \le a < 2\alpha, b< 1/a.$ In particular, any B-spline $B_N, N\ge 2,$ generates a frame $\{E_{mb}T_{na}B_N\}_{m,n\in \mz}$ whenever $N/2 \le a < N, \, b< 1/a.$ This confirms a recent conjecture by Gr\"ochenig in that particular region. Inspired by this result we prove that the full conjecture  holds if it can be verified in the region determined by the inequalities
$1/2 \le ab <1, \, a<N/2.$

The key result in the paper is Theorem \ref{2-2-1}, which
characterizes the frame property of $\mts$ in the aforementioned
region. The proof is quite complicated and is split into several
lemmas and intermediate steps.  The idea of the proof was gained
through the work on the special case with translation parameter
$a=1$ (see the paper \cite{CKK10}), as well as the observation
that the duality condition \eqref{14-20} forces a certain behavior
of the window  $g$ around points $x_0+a$ for which $g(x_0)=0.$ As
further help to understand the idea behind the proof we prove the
steps directly in a concrete case, see Example \ref{4108e}. For
more informations about Gabor systems and frames we refer to the
monographs [5,1].

\section{General results}
Given $ \al>0$, let \bee \notag V_\alpha:= \{f\in C(\mr) \ | \
\supp \, f=[-\alpha,\alpha],
\ \mbox{$f$ has a finite number of zeros on} \ [-\alpha,\alpha]\}. \\
\ \label{wi-1} \ \ene

We first characterize the frame property of  $\mts$ for
functions $g \in V_\al$ and points $(a,b)$ in the region in $\mr_+^2$
determined by the inequalities $\alpha \le a <2\alpha, b< 1/a.$ In
order to do this, we need to introduce some parameters and other
tools. Consider $(a,b)$ belonging to the described region, and choose $M\in
\mn$ such that $ab\in [\frac{M-1}{M},\frac{M}{M+1}[.$
Let $\kappa$ be the largest integer for which
  $(1-ab)\kappa \leq b \al$.
Then $0\leq \kappa \leq M-1$ because
$$\kappa \leq \frac{b\al}{1-ab} \leq \frac{ab}{1-ab}<\frac{M}{M+1}\left(1- \frac{M}{M+1}\right)^{-1}=M.$$

If $\kappa \neq 0$, let $n\in\{1,2,\cdots, \kappa\},$ and define the function
$R_{n}$ on (a subset of) $]a-\al,\al-(1-ab)\frac{n}{b} ]
\subset ]a-\al, \al]$ by \bee \label{4277a} R_{n}(y):=
\frac1{g(y)} \prod_{k=1}^{n-1} \frac{
g(y+(1-ab)\frac{k}{b}-a)}{g(y+(1-ab)\frac{k}{b})} , \ \
n=1,2,\cdots, \kappa. \ene We use the standard convention that the
empty product is $1$. It is easy to see that $R_{n}$ indeed is defined on $]a-\al,
\al-(1-ab)\frac{n}{b}],$ except possibly on  a finite set of
points.
Similarly, still if $\kappa \neq 0$, for $n\in\{1,\cdots, \kappa\}$ we define the function
$L_{n}(y)$ on   (a subset of) $[-\al + (1-ab)\frac{n}{b},\al-a[
\subset [-\al, \al-a[$ by
$$L_{n}(y):=
 \frac1{g(y)}\prod_{k=1}^{n-1} \frac{  g(y-(1-ab)\frac{k}{b}+a )}{
 g(y -(1-ab)\frac{k}{b})}, \ \
n=1,2,\cdots, \kappa.
$$

We now state the announced characterization of the frame property.

\bt \label{2-2-1} Let $g\in V_\alpha$ for some $\alpha>0$ and
assume that $ \al\leq a < 2 \al$ and
$ab\in[\frac{M-1}{M},\frac{M}{M+1}[$ for some
$M\in\mathbb{N}\setminus\{1\}.$ Let $\kappa\in \{0,1,\cdots,
M-1\}$ be the largest integer for which    $(1-ab)\kappa \leq
b\alpha$.
Then $\{E_{mb} T_{na} g \}_{m,n\in \mathbb{Z}}$ is a Gabor frame if and only if the following  conditions are satisfied:
\begin{itemize}
  \item[{\rm(i)}] $|g(x)|+|g(x+a)|>0, \ x\in [-a,0]$;
    \item[{\rm(ii)}] If $\kappa\neq 0$  and if  there exist $n_+\in\{1,2,\cdots, \kappa\}$ and
                    $y_+\in ]a-\al, \al- (1-ab)\frac{n_+}{b}]$ such that
                    $g(y_+)=0$ and
               $ \lim_{y\rightarrow y_+} |R_{n_+}(y)|=\infty$,
               then
               \begin{equation*}
                g(y_++ (1-ab)\frac{n_+}{b}-a)\neq 0;
               \end{equation*}
    \item[{\rm(iii)}]
If $\kappa\neq 0$ and if there exist $n_-\in\{1,2,\cdots, \kappa\}$
and
                    $y_-\in[-\al +  (1-ab)\frac{n_-}{b},\al-a[$  such that
                    $g(y_-)=0$ and  $ \lim_{y\rightarrow y_-} |L_{n_-}(y)|=\infty$,
                    then \bes g(y_--(1-ab)\frac{n_-}{b}+a)\neq 0;\ens
\item[{\rm(iv)}] For $y_+,y_-, n_+, n_-$ as in {\rm(ii)} and
{\rm(iii)},
$$ y_++(1-ab)\frac{n_+}{b} \neq y_--(1-ab)\frac{n_-}{b}+a. $$
\end{itemize}
In the affirmative case, there exists a  dual window $h$
    with {\em $\supp\ h\subseteq [-aM,aM]$}.
\et

We remark that if $\kappa=0$ then the conditions (ii)-(iv) are trivially satisfied. We also note that Theorem \ref{2-2-1} is similar, but significantly
more general than Theorem 2.3 in \cite{CKK10}. The main difference is that in the current paper  the support size of $g$
(measured by the parameter $\alpha$) and
the translation parameter $a$ can vary, subject to the
restriction $ \al\leq a < 2 \al;$  on the other hand
\cite{CKK10} dealt with the case $\alpha=a=1.$ This
modification turns out to be
instrumental for our applications to B-splines.

The proof of the necessity of the conditions in
Theorem \ref{2-2-1} is similar to the proof in
\cite{CKK10}, so we skip this part.  On the other hand,
it requires much more work to prove that
$\mts$ is a frame if the conditions (i)-(iv)
are satisfied. We prove this part of the theorem in the appendix. In the subsequent example
we prove directly that a certain Gabor system is a frame, following the steps from the
proof of the general result; the hope is that the analysis of this concrete case will help the reader to understand the idea behind the general proof.


\bex \label{4108e} Let $\al=9/10$ and consider a function $ g\in
V_\al,$ having the single zero $1/5$ within $]-1,1[.$ Let $a=1$
and $b=3/5$. Then  $ab\in[\frac{M-1}{M},\frac{M}{M+1}[$ for $M=2$.
Note that $ (1-ab)/b = 5/3 -1 = 2/3  \leq 9/10 =\al.$ This implies
that  $\kappa= 1.$ Trivially, $|g(x)|+|g(x+a)|>0, \ x\in [-a,0]$.
Let $n_+:=1$  and $y_+:=1/5$. Then  $y_+\in  ]a-\al, \al -
(1-ab)\frac{n_+}{b}]=]1/10, 7/30]$ and $g(y_+)=0.$  Furthermore,
$R_{n_+}(y)= g(y)^{-1},$ so
               $ \lim_{y\rightarrow y_+} |R_{n_+}(y)|=\infty.$ It is also clear that
$g(y_+ +(1-ab)\frac{n_+}{b}-a)=g(-2/15)  \neq 0$.

It is an easy consequence of the duality conditions for Gabor
frames (see \eqref{14-20} in the Appendix) that two real valued,
bounded functions  $g, h\in L^2(\mr)$ with $\supp \ h \ \subseteq
[-aM,aM]=[-2,2]$ generate dual frames $\{E_{mb}T_{na}g\}_{m,n\in \mz}$
and $\{E_{mb}T_{na}h\}_{m,n\in \mz}$ if and only if  for
 $n=0,\pm 1$ and
 $a.e.$ $x\in [\frac{n}{b}-a,\frac{n}{b}]$,
\begin{equation}\label{14-16}
 g(x-\frac{n}{b})h(x) +g(x-\frac{n}{b}+a)h(x+a)=b \delta_{n,0}.
\end{equation}
We will check \eqref{14-16} directly following the steps in the
general proof of Theorem \ref{2-2-1}. Motivated by a general result, see Lemma \ref{b-1},
we choose to put $h(x)=0$ for  $x\notin [-a-\al,-\frac{1}{b}] \cup
[-\al,\al] \cup [\frac{1}{b}, \al+a].$
Then $h(x)=h(x+a)=0$ for
  $x\in ]\al, \frac{1}{b}[$,     which is a subinterval of
$[\frac{1}{b}-a,\frac{1}{b}]$; thus \eqref{14-16} holds for $n=1$
and  $x\in ] \al, \frac{1}{b}[$. Similarly, \eqref{14-16} holds for
$n=-1$ and  $x\in] -\frac{1}{b}-a, -a -\al[$.
Note that $g(x-\frac{1}{b}+a)=0$ if and only if
  $x=y_+ + \frac{1}{b}-a $.
Let us for a moment assume that $h$ is
chosen on $[\frac{1}{b}-a, \al]$ as a bounded function such that
 $h$ is continuous at $y=y_++\frac{1}{b}-a$ and
\begin{equation} \label{14-17}
\lim_{y\rightarrow y_+}
 \left\{ h(y+(1-ab)\frac{n_+}{b})R_{n_+}(y)
 \right\}
\end{equation}
exists; letting $x=y + \frac{1}{b}-a$, this means that
\begin{equation*}
\lim_{x\rightarrow y_+ +\frac{1}{b}-a}
 \left\{\frac{h(x)}{g(x-\frac{1}{b}+a)}
 \right\}
\end{equation*}
exists.
Then, defining $h$ on $[\frac{1}{b}, a+\al]$ by
$$h(x+a)=\left\{
\begin{array}{ll}
-\dfrac{g(x-\frac{1}{b})h(x)}{g(x-\frac{1}{b}+a)}, & x \in [\frac{1}{b}-a, \al]\setminus \{ y_+ +\frac{1}{b}-a\}; \\
-\lim_{t\rightarrow y_+ +\frac{1}{b}-a} \left\{
\dfrac{h(t)}{g(t-\frac{1}{b}+a)}
 \right\}     g(x-\frac{1}{b}), & x = y_+ +\frac{1}{b}-a,
\end{array}
\right.
$$
\eqref{14-16} holds for $n=1$ and $x\in [\frac{1}{b}-a, \al]$.
Hence
$\mts$ is a frame if we  can define $h$ as  a bounded function on
$[-\alpha, \alpha]$ such that
\begin{itemize}
    \item[(a)]  $h$ is continuous at $y=y_+ + \frac{1}{b}-a $ and
    \eqref{14-17}  holds;
    \item[(b)] the duality condition \eqref{14-16} holds for $n=0$ and $x\in [-a,0],$ $i.e.$,
        \bee \label{4925g} g(x)h(x)+ g(x+a)h(x+a)=b, \, x\in [-a,0].\ene
\end{itemize}
  Let $\tilde y_+:=y_+ + (1-ab)\frac{1}{b}$ and let $B_1:= ]y_+
-\epsilon, y_+ +\epsilon[ \cup ]\tilde y_+  -\epsilon, \tilde y_+
+\epsilon[ \cup ]\al-\epsilon,\al+\epsilon[$ and
$B_2:=]-\al-\epsilon,-\al+\epsilon[,$
for an $\epsilon >0$ chosen such that
\bei \item[(i)] $|g(x)|\ge \delta > 0$ for $x\in (B_1-a)\cup (B_2+a)$ and
some $\delta>0$;
 \item[(ii)] $B_1\cap (B_2 +a)=\emptyset.$
\eni
Note that $g(x)\neq 0, \ x\in [\al-a,a-\al].$
By  continuity of $g$, $\inf_{x\in [\al-a,a-\al]} |g(x)|>0$.  We
define $h(x):=\frac{b}{g(x)}, \ x\in [\al-a,a-\al]$, which is thus
a bounded function. Note that for $x\in [-a,-\al]$, we have
$g(x)=0$, and therefore
\begin{equation}\label{14-18}
g(x)h(x)+g(x+a)h(x+a)=b.
\end{equation}
Similarly, \eqref{14-18} holds for $x\in [\al-a,0]$, $i.e.$, we have
now verified (b) on the subinterval $[-a,-\al]\cup[\al-a,0]$. We
now put $h=0$ on $B_1\cap[a-\al,\al]$. Then
$h$ is continuous at $y=y_+ + \frac{1}{b} -a $  and
 $ \lim_{y\rightarrow
y_+} \left\{ h(y     +(1-ab)\frac{1}{b} )R_{1}(y) \right\}=0. $
Hence  (a)  holds. We define $h$ on $( B_1 -a )\cap [-\al,\al-a]$
by $ h(x)= \frac{b-g(x+a)h(x+a)}{g(x)}= \frac{b}{g(x)};$ thus $h$
is bounded here by the choice of $\epsilon$,  and (b) holds on $(
B_1 -a )\cap [-\al,\al-a]$. Similarly, put $h=0$ on $B_2\cap
[-\al,\al-a]$ and define $h$ on $(B_2+a)\cap [a-\al,\al]$ by $
h(x)= \frac{b}{g(x)};$ thus $h$ is bounded,  and (b) holds on $B_2
\cap [-\al,\al-a]$.
We finally put $h=0$ on $[-\al, \al-a]\setminus ((B_1-a)\cup
B_2).$ Note that the zeroset of $g$ within $[-\al,\al]$ is
$\{-\al,y_+,\al\}$, so $g(x)\neq 0$ for $x\in
\overline{[a-\al,\al]\setminus (B_1\cup(B_2+a))};$ using the
continuity of $g$ implies that $\inf_{x\in [a-\al,\al]\setminus
(B_1\cup(B_2+a))} |g(x)|>0$. We define $h(x)=\frac{b}{g(x)},\ x\in
[a-\al,\al]\setminus (B_1\cup(B_2+a));$  thus, we have now defined
$h$ everywhere as a bounded function, and (b) holds for $x\in
[-\al, \al-a]\setminus ((B_1-a)\cup B_2).$    This completes the
proof of (b), and hence the proof of $\mts$ being a Gabor frame
with a dual window supported on $[-2,2].$
 \ep \enx


From Theorem \ref{2-2-1} we can immediately extract  the possible obstruction curves, i.e., the points $(a,b)$ for which a given function $g\in V_\alpha$ might not generate a frame $\mts.$
Assume that $g\in V_\alpha$ satisfies the standing assumptions  in Theorem \ref{2-2-1} as well as the condition
\begin{equation}\label{14-10}
|g(x)|+|g(x+a)|>0,\ x\in [-a,0].
\end{equation}
Then, if $\kappa \neq 0,$ the possible obstructions take place on the curves determined by the equations

\begin{eqnarray}
&&  y_++(1-ab)\frac{n_+}{b}-a=y_- , \label{14-11}\\
&& y_--(1-ab)\frac{n_-}{b}+a=y_+ , \label{14-12}\\
&&  y_++(1-ab)\frac{n_+}{b}= y_--(1-ab)\frac{n_-}{b}+a.
\label{14-13}
\end{eqnarray}for some $y_+, y_-, n_+, n_-$ as in the theorem.
The equations \eqref{14-11} and \eqref{14-12} both take the form
\bee \label{4108a} b= \frac{n}{y_--y_++an+a}\ene for some $n\in\{1,2,\cdots, \kappa\},$ while
\eqref{14-13} means that
\bee \label{4108b} b= \frac{n_-+n_+}{y_--y_++(n_-+n_+)a+a}\ene for some
$n_-, n_+\in\{1,2,\cdots, \kappa\}.$
Note that these curves only depend
on the location of the zeros of the function $g\in V_\alpha,$ not on the specific function.

Interestingly, the equations \eqref{4108a} and \eqref{4108b} show that for functions
$g\in V_\alpha$ the obstructions take place on ``hyperbolic curves:" this is
similar to the result in \cite{LyNes}, where Lyubarski and Nes  showed that for any odd function in the Feichtinger algebra $M^1,$ (in particular,
the first order Hermite function) the
points $(a,b)$ for which $ab=1-1/M= \frac{M-1}{M}$ for $M=2,3,\dots$ do not belong to the frame set.

For functions $g\in V_\alpha$ with no zeroes in $]-\alpha, \alpha[$ the conditions in
Theorem \ref{2-2-1} are clearly satisfied, which yields the following:

\bc \label{4714a}
Let $\alpha >0$. Assume that $g$ is a continuous  function with {\em $\supp$}
$g = [-\al, \alpha]$, and that
$$ g(x)>0, \ \ x\in ]-\al, \alpha[.$$
Then   $\{E_{m b}T_{n a}  g\}_{m,n\in \mz}$ is a frame whenever
$\alpha \le  a < 2\alpha, \, 0<b< 1/a.$
\ec

\section{B-splines and a conjecture by Gr\"ochenig} \label{50131a}

Let us now consider the B-splines $B_N, \, N\in \mn,$ defined recursively by
\bes B_1 = \chi_{[-1/2, 1/2]}, \, B_{N+1} = B_N *B_1.\ens  The frame properties of
$\{E_{mb}T_{na}B_1\}_{m,n\in \mz}$ are well known (see the work by Janssen \cite{Jan7}
and \cite{Sun} by Dai and Sun which finally solved the
so-called $abc$-problem), so we focus on the case $N\ge 2,$ where $B_N$ is a continuous function supported on $[-N/2, N/2].$ Furthermore
the function $B_N, \, N\ge 2,$ is strictly positive on the interval  $]-N/2, N/2[,$ so
Corollary \ref{4714a} implies that  $\{E_{mb}T_{na}B_N\}_{m,n\in \mz}$
is a frame whenever $N/2 \le  a < N, \, 0<b< 1/a.$  Several other  results about the frame set ${\cal F}(B_N)$ are known. We collect them here for easy reference:

\bpr \label{new} Let $N\in \mn \setminus \{1\},$ and consider $a,b>0$ such that
$ab<1.$ Then the following hold:

\bei \item[{\rm (i)}] $\{E_{mb}T_{na}B_N\}_{m,n\in \mz}$ is not a frame if $a\ge N.$
\item[{\rm (ii)}] $\{E_{mb}T_{na}B_N\}_{m,n\in \mz}$ is not a frame if $b=2,3,\dots .$
\item[{\rm (iii)}] $\{E_{mb}T_{na}B_N\}_{m,n\in \mz}$ is a frame if $a< N, \, b\le 1/N.$
\item[{\rm (iv)}] $\{E_{mb}T_{na}B_N\}_{m,n\in \mz}$ is a frame if there exists
$k\in \mn$ such that \bee \label{4710d}   1/N < b < 2/N, \,  N/2 \leq  ak < 1/b.\ene
\item[{\rm (v)}] $\{E_{mb}T_{na}B_N\}_{m,n\in \mz}$ is a frame if $b\in \{1, \frac12, \dots, \frac1{N-1}\}.$
\item[{\rm (vi)}] $\{E_{mb}T_{na}B_N\}_{m,n\in \mz}$ is a frame if $a=\frac{k}{p}
$ for some $k=1, \dots, N-1, \, p\in \mn,$ and $b< 1/k.$
\eni \epr

\bp The results in (i) and (iii) are classical.
Also (ii) is a well known result, originally due to Del Prete \cite{Del1} and rediscovered in
\cite{GJ}. For $k=1,$
the statement in (iv) is a consequence  of Corollary \ref{4714a}.
In general, if \eqref{4710d} holds for some $k\in \mn \setminus \{1\},$ then
this implies that $\{E_{mb}T_{nka}B_N\}_{m,n\in \mz}$ is a frame, and we infer that the larger system
$\{E_{mb}T_{na}B_N\}_{m,n\in \mz}$ itself is a frame
(because the upper bound holds automatically).
The result in (v) was recently proved by  Kloos and St\"ockler  \cite{KlSt}, who also proved (vi) for $p=1;$ The case of $p\in \mn$ in (vi) yields an oversampling
of the case $p=1,$ and therefore a frame.
\ep

Based on (i)--(iii) in Proposition \ref{new} Gr\"ochenig formulated a conjecture
about the frame set ${\cal F}(B_N)$ in \cite{G8}. Basically it says that the frame set
consists of all the points $(a,b)\in \mr_+^2$ that avoids the known obstructions:

\vn {\bf Conjecture} For any $N\ge 2,$
\bes {\cal F}(B_N)= \{(a,b)\in \mr_+^2 \, \big| \,  a<N, \, ab <1, b \neq 2, 3, \dots\}.\ens

We will now show that the  conjecture is true if
we can prove the frame property in a certain `` hyperbolic strip."

\begin{figure}
\includegraphics[width=5in,height=2.4in]{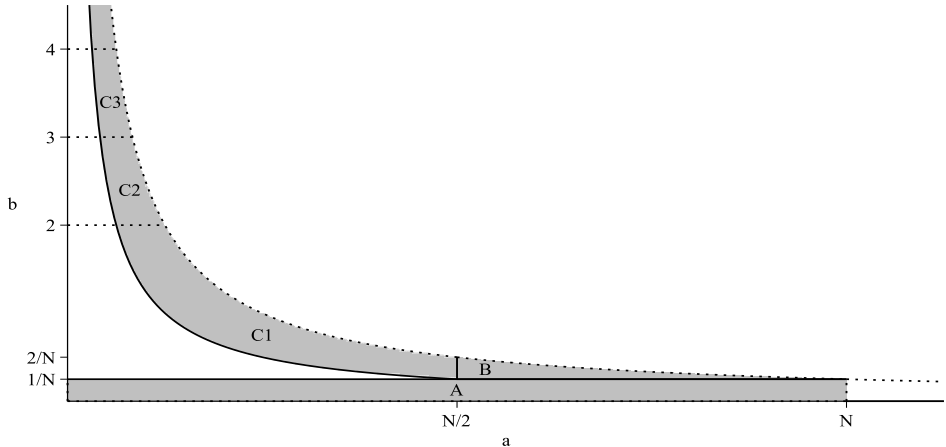}\hfil
\caption{The set $A$ belongs to the frame
set for $B_N, N>1.$ Corollary \ref{4714a} proves that $B$ also
belongs to the frame set (see also the introduction to Section \ref{50131a}); and by Proposition \ref{4710b}
the conjecture by Gr\"ochenig is true if it can be verified in the
regions $C1, C2, \dots.$} \label{4108f}
\end{figure}

\bpr \label{4710b} The conjecture is true if  $\{E_{mb}T_{na}B_N\}_{m,n\in \mz}$ is a frame for all $(a,b)\in \mr_+^2$
for which
\bee \label{4710af} a<N/2, \, 1/2 \le ab <1, \, b\notin \{2,3,\dots\}.\ene
 \epr

\bp To get a geometric understanding we refer to Figure \ref{4108f}. We note that
Corollary \ref{4714a} confirms the frame property in the region determined by
the inequalities $N/2<a <N, ab<1;$ furthermore the frame property is satisfied
for $a<N, b\leq 1/N$  (i.e., the region A on Figure \ref{4108f}). Thus,
it suffices to show that the parameter region
determined by the inequalities
\bes 0 < ab < \frac12, \, \, \, \frac1{N} < b \notin \{2, 3, \dots\}\ens
is contained in the frame set ${\cal F}(B_N)$ under the given assumption. Note that $\{E_{mb}T_{na}B_N\}_{m,n\in \mz}$ is a Bessel sequence for all $a,b>0,$ i.e.,  we only need to check the lower frame condition.

Since  $0 < 2ab <1,$ choose the unique $M\in \mn$ such that $\frac1{M+1} \le 2ab< \frac1{M}.$ By splitting into the cases $2Ma<N/2$ and $2Ma\ge N/2$ it follows
that the system $\{E_{mb}T_{n2Ma}B_N\}_{m,n\in \mz}$ is a frame; this clearly
implies that $\{E_{mb}T_{na}B_N\}_{m,n\in \mz}$  satisfies the lower frame bound as well. \ep

\section*{Appendix: Proof of Theorem \ref{2-2-1}} \label{s5-1}
Let $M\in \mn$, and assume that $\frac{M-1}{M}\leq
ab<\frac{M}{M+1}.$ The starting point is the duality conditions by
Ron \& Shen \cite{RoSh, RoSh5} and Janssen \cite{J}, which by an
easy calculation implies that two real valued, bounded functions
$g,h\in \ltr$ with $\supp\ g \subseteq [-a,a]$, $\supp\ h
\subseteq [-aM,aM]$, generate dual frames $\{E_{mb}T_{na}g\}_{m,n\in
\mz}$ and $\{E_{mb}T_{na}h\}_{m,n\in \mz}$ if and only if for
 $n=0,\pm 1,\pm 2, \cdots, \pm(M-1)$ and
$a.e.$  $x\in [\frac{n}{b}-a,\frac{n}{b}]$,
\begin{equation}\label{14-20}
 g(x-\frac{n}{b})h(x) +g(x-\frac{n}{b}+a)h(x+a)=b \delta_{n,0}.
\end{equation}

We will now consider a function $g\in V_\alpha$ that satisfies the
conditions (i)--(iv) in Theorem \ref{2-2-1}. We will prove that
$\mts$ is a frame by constructing a dual window $h.$
In the following lemma, we use the insight gained from the proofs in \cite{CKK10}
to define $h$ on certain intervals, in such a way that
\eqref{14-20} is satisfied for some of the
relevant values of $n$ and on certain intervals.  After that the subsequent lemma states conditions that yields a definition of $h$ on the remaining parts of the real line in such a way that all the duality conditions are satisfied.

\begin{figure}
\includegraphics[width=5.5in,height=1in]{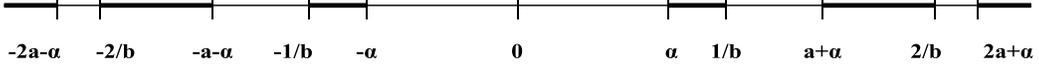}\hfil
\caption{The figure shows the set where $h$ is defined
to vanish by \eqref{b-2}, in the case $\kappa=2$.  }
\end{figure}

\begin{lemma}\label{b-1}
Let $\al, a,b >0$ be given such that $  \al\leq a <2\al$ and
$ab\in [\frac{M-1}{M},\frac{M}{M+1}[$ for some $M\in\mathbb{N}\setminus\{1\}$.
Assume  that
\begin{equation}\label{b-2}
   h(x)=0, \ \ x\notin -\left(\bigcup_{k=1}^{\kappa}[\frac{k}{b},ak+\alpha]\right)
   \cup [-\al,\al]
   \cup \bigcup_{k=1}^{\kappa}[\frac{k}{b},ak+\alpha].
\end{equation}
Then the following hold:
\begin{itemize}
    \item[{\rm (a)}]
$h(x)=h(x+a)=0$ for $n=1,\cdots,\kappa$ and $x\in
]\alpha+a(n-1),\frac{n}{b}[$,
and
for $n=-1,\cdots,-\kappa$ and
$x\in ]\frac{n}{b}-a,an-\al[$
;
    \item[{\rm (b)}]
$h(x)=h(x+a)=0$  for $n=\pm(\kappa+1), \cdots, \pm(M-1)$
and $x\in [\frac{n}{b}-a,\frac{n}{b}].$
\end{itemize}
\end{lemma}
\bp
Note that for
$n=1,2,\dots, M-1,$ $ab\ge \frac{n}{n+1};$ thus,
$\frac{n}{b}-a \leq a n <\frac{n}{b}.$

\noindent (a): For $1\leq n \leq \kappa$, we note that the statement in (a) only
involves the function values of $h$ for $x\in ]\al + a(n-1),n/b[$ and $x+a \in ]\al+an, n/b+a[.$
Since $\left( ]\al + a(n-1),n/b[\cup ]\al+an, n/b+a[  \right)\cap \supp\ h=\emptyset$,
(a) holds for $1\leq n\leq \kappa$.
Similarly,
(a) holds for $-\kappa \leq  n\leq -1$.

 \noindent (b): For $\kappa+1 \leq  n \leq M-1$, the statement in (b)  only
involves the values of $h$ for $x\in [\frac{n}{b}-a,\frac{n}{b}]$ and $x+a \in [\frac{n}{b},\frac{n}{b}+a]$;
by the definition of $\kappa$,
    we have $b\alpha< (\kappa +1)(1-ab)$, $i.e.$,
$a\kappa+\alpha < \frac{\kappa+1}{b}-a$; thus 
$\left( [\frac{n}{b}-a, \frac{n}{b}+a ] \right)\cap \supp \ h =\emptyset.$
Hence  \eqref{14-20} holds for $\kappa+1\leq n\leq M-1.$
Similarly, (b) holds for $-M+1\leq n \leq -\kappa-1$.   \ep

Note that condition (b) in Lemma \ref{14-20} is empty if $\kappa=0$; condition (c) is empty if $\kappa=M-1$.

 By Lemma \ref{b-1}, we see that
\eqref{14-20} holds for $n=1,\cdots,\kappa$ and $x\in
]\alpha+a(n-1),\frac{n}{b}[$,
and
for $n=-1,\cdots,-\kappa$ and
$x\in ]\frac{n}{b}-a,an-\al[$. Similarly,
\eqref{14-20} holds  for $n=\pm(\kappa+1), \cdots, \pm(M-1)$
and $x\in [\frac{n}{b}-a,\frac{n}{b}].$
What remains is to show that we can define $h$ on the set
\bes -\left(\bigcup_{k=1}^{\kappa}[\frac{k}{b},ak+\alpha]\right)
   \cup [-\al,\al]
   \cup \bigcup_{k=1}^{\kappa}[\frac{k}{b},ak+\alpha]\ens such that \eqref{14-20} holds for
$n=1, \cdots,\kappa$ and $x\in [\frac{n}{b}-a,\al+a(n-1)]$, and for
$n=-1, \cdots,-\kappa$ and $x\in [an-\al, \frac{n}{b}]$,
as well as for $n=0$ and $x\in[-a,0].$ The following lemma states sufficient conditions for the first of these requirements to be satisfied.
The result is a   minor adaption of Lemma 3.3 in \cite{CKK10}, so the proof is omitted.

\begin{lemma}\label{b-3}
Let $\alpha, a,b>0$ be given such  that  $\alpha \le a < 2\alpha$
and $ab\in [\frac{M-1}{M},\frac{M}{M+1}[$ for some $M\in\mathbb{N}\setminus\{1\}$. Let $g\in V_\alpha$, and assume  that
$g(x) \neq 0$ for $x\in [\alpha -a, a - \alpha].$
 Assume further that $\kappa \neq 0$ and
that $h$ is chosen
 on  $[-\al,\al]$ as a bounded function  such that the following
conditions hold:
\begin{itemize}
    \item[\rm{(1)}] If there exist $n_+\in\{1,2,\cdots, \kappa\}$ and $y_+\in ]a-\al,  \al-(1-ab)\frac{n_+}{b}]$
such that $g(y_+)=0$ and $\lim_{y\rightarrow y_+}
|R_{n_+}(y)|=\infty $, then
 $h$ is continuous at $y=y_+ + (1-ab)\frac{n_+}{b}$  and
the limit
\begin{equation} \label{b-4}
\lim_{y\rightarrow y_+}
 \left\{ h(y+(1-ab)\frac{n_+}{b})R_{n_+}(y)
 \right\}
\end{equation}
exists;
    \item[\rm{(2)}] If there exist $n_-\in\{1,2,\cdots, \kappa\}$ and
    $y_-\in [-\al+ (1-ab)\frac{n_-}{b},\al-a]$
such that $g(y_-)=0$ and $\lim_{y\rightarrow y_-}
|L_{n_-}(y)|=\infty $, then
 $h$ is continuous at $y=y_- - (1-ab)\frac{n_-}{b}$  and
the limit
\begin{equation} \label{b-5}
\lim_{y\rightarrow y_-} \left\{ h(y-(1-ab)\frac{n_-}{b})
L_{n_-}(y)\right\}
\end{equation}
 exists.
\end{itemize}
Then the function $h$ can be defined on the interval
 $  -(\bigcup_{k=1}^{\kappa}[\frac{k}{b},ak+\alpha] ) \cup
 \bigcup_{k=1}^{\kappa}[\frac{k}{b},ak+\alpha]$ such that
the duality condition \eqref{14-20} holds  for
$n=1,\cdots, \kappa$ and $ x\in
[\frac{n}{b}-a, \al+a(n-1)],$
as well as for
$n=-1,\cdots,-\kappa$ and $ x\in
[an-\al, \frac{n}{b}];$
the function $h$ is  bounded, and the values of $h$ outside $ -(\bigcup_{k=1}^{\kappa}[\frac{k}{b},ak+\alpha] ) \cup
[-\al,\al]\cup \bigcup_{k=1}^{\kappa}[\frac{k}{b},ak+\alpha]  $ are irrelevant.
\end{lemma}

\begin{rem}{\em
In Lemma \ref{b-3}, if $y_+$ is the end point of the interval
$]a-\al,  \al-(1-ab)\frac{n_+}{b}]$, $i.e.$,
$y_+=\al-(1-ab)\frac{n_+}{b}$,  the limit $\lim_{y\rightarrow
y_+}$ in \eqref{b-4} should be understood as the limit from the
left; similarly, if $y_-=-\al+(1-ab)\frac{n_-}{b}$, the limit
$\lim_{y\rightarrow y_+}$ in \eqref{b-5} should be understood as
the limit from the right.
\ep}
\end{rem}

We can now  complete the proof of the sufficiency in
Theorem \ref{2-2-1}:

\noindent{\bf Proof of   Theorem
\ref{2-2-1}:} Assume that the conditions (i)--(iv) in Theorem \ref{2-2-1} hold.
Note that $g(x)=0,\ x\in [-a,-\al]\cup[\al,a]$, since $g\in V_\al$.
This together with  condition (i) in Theorem \ref{2-2-1} implies that
\begin{equation}\label{14-14}
  g(x)\neq 0, \ x\in [\al-a,a-\al].
\end{equation}
Following  \eqref{b-2}, let
\bes
   h(x):=0, \ \ x\notin -\left(\bigcup_{k=1}^{\kappa}[\frac{k}{b},ak+\alpha]\right)
   \cup [-\al,\al]
   \cup \bigcup_{k=1}^{\kappa}[\frac{k}{b},ak+\alpha].
\ens
Via Lemma \ref{b-3} and the
comment just before the lemma,  $\mts$ is a frame
if we  can define
$h$ as a bounded function on $[-\alpha, \alpha]$ in such a way that
\begin{itemize}
    \item[(a)] the conditions in Lemma \ref{b-3} (1) and (2) hold;
    \item[(b)] the duality condition \eqref{14-20} holds for $n=0$ and $x\in [-a,0],$ $i.e.$,
        \bee \label{4925g} g(x)h(x)+ g(x+a)h(x+a)=b, \, x\in [-a,0].\ene
\end{itemize}
We will split the definition of $h$ on $[-\al,\al]$ into several intervals. In fact, we will first
define $h$ on $[\al-a,a-\al]$ and then on small balls around certain shifts of the zeros. First, we need
some notation. For $m,n=0,1,\cdots, \kappa$,
we define the sets $Y_n$ and $W_m$  by
\begin{eqnarray*}
&&Y_0= \{y_{0,i}\in ]a-\al,\al ] \ : g(y_{0,i})=0\}_{i=1,2,\cdots,r_0} \\
&&Y_n=\{y_{n,i}\in ]a-\al,\al-(1-ab)\frac{n}{b} ] \ : g(y_{n,i})=0
\text{ and $\lim_{y\rightarrow y_{n,i}} |R_{n}(y)|=\infty
$}\}_{i=1,2,\cdots,r_n}
\end{eqnarray*}
and
\begin{eqnarray*}
&&W_0=\{w_{0,j}\in [-\al,\al-a [ \ : g(w_{0,j})=0\}_{j=1,2,\cdots,l_0}\\
&&W_m=\{w_{m,j}\in [-\al+(1-ab)\frac{m}{b},\al-a [ \ : g(w_{m,j})=0 \text{ and
$  \lim_{y\rightarrow w_{m,j}} |L_{m}(y)|=\infty$}\}_{j=1,2,\cdots,l_m}
\end{eqnarray*}
where $r_n$ and $l_m$ are the cardinalities of $Y_n$ and $W_m$,
respectively. In words: since $g(x)\neq 0$ for $x\in [\al-a,a-\al]$,
 the sets $Y_0$ and $W_0$ yield enumerations of the zeros for $g$
 within $[-\al,\al],$
split into the positive, respectively, negative part;
the sets $Y_n$ and $W_n, n\ge 1,$ yield
enumerations of selected zeros within certain subsets of $[-\al,\al].$

We denote the open interval of radius $r>0$
centered at $x$ by $B(x;r)= ] x- r, x+r[.$
For $y_{n,i} \in Y_n$, $w_{m,j} \in W_m$ for $n,m=0,1,\cdots, \kappa$,
let
$\tilde y_{n,i}:=y_{n,i}+(1-ab)\frac{n}{b},\
\hat w_{m,j}:=w_{m,j}-(1-ab)\frac{m}{b}.$
If $n,m \geq 1$, then  by the conditions (ii), (iii) and (iv)
 in Theorem \ref{2-2-1}(3), we have
\bee \label{4918a} g(\tilde y_{n,i}-a)\neq 0 \neq g(\hat w_{m,j}+a),
\, \, \mbox{and} \, \,
\tilde y_{n,i} \neq \hat w_{m,j}+a.\ene
Note that $g(\tilde y_{0,i})=g(\hat w_{0,j})=0$. Then we also have
$\tilde y_{0,i} \neq \hat w_{m,j}+a,\ \tilde y_{n,i}-a \neq \hat w_{0,j}$ for  $m,n\geq 1$, and
$g(\tilde y_{0,i}-a)\neq 0 \neq g(\hat w_{0,j}+a)$ by the condition (i) in Theorem \ref{2-2-1}; thus, \eqref{4918a} actually holds for all
$m,n=0,1,\cdots, \kappa.$
Then we can choose $\epsilon >0$ so that
\bei \item[(i)] $|g(x)|\ge \delta > 0$
for  $x\in
 B(\tilde y_{n,i}-a; \epsilon)
\cup B(\hat w_{m,j}+a ;\epsilon)
$ and some $\delta>0$;
\item[(ii)] For $m,n=0,1,\cdots, \kappa$,  and
$i=1,2,\cdots,r_n,
j=1,2,\cdots,l_m,$
\begin{equation}\label{epsilon2-1}
B(\tilde y_{n,i}; \epsilon) \cap B(\hat w_{m,j}+a ;\epsilon)=\emptyset.
\end{equation}
\eni


\noindent {\bf Definition of $h$ on $[\al-a,a-\al]$:}
By  \eqref{14-14} and continuity of $g$,
$\inf_{x\in [\al-a,a-\al]} |g(x)|>0$.  We define
$h(x):=\frac{b}{g(x)}, \ x\in [\al-a,a-\al]$, which is thus a
bounded function. Note that for $x\in [-a,-\al]$, we have
$g(x)=0$, and therefore
\begin{equation}\label{b-6}
g(x)h(x)+g(x+a)h(x+a)=b.
\end{equation}
Similarly, \eqref{b-6} holds for $x\in [\al-a,0]$, $i.e.$, we have now verified (b)
on the subinterval $[-a,-\al]\cup[\al-a,0]$.


\noindent {\bf Definition of $h$ on $B(\tilde y_{n,i};\epsilon)\cap[a-\al,\al]$:}
On this interval, put $h=0$.
If
$1\leq n\leq \kappa$, then,
$h$ is continuous at $y=\tilde y_{n,i}$  and
\begin{equation}\label{b-9}
 \lim_{y\rightarrow y_{n,i}}
\left\{ h(y     +(1-ab)\frac{n}{b} )R_{n}(y) \right\}=0.
\end{equation}
 Hence the condition in Lemma \ref{b-3} (1)  holds.


\noindent {\bf Definition of $h$ on $B(\tilde y_{n,i}-a;\epsilon)\cap[-\al,\al-a]$:}
 We define $h$ on this set
by $ h(x)= \frac{b-g(x+a)h(x+a)}{g(x)}= \frac{b}{g(x)};$ thus $h$ is bounded
here by the choice of $\epsilon$,  and (b) holds on
$B(\tilde y_{n,i}-a;\epsilon)\cap[-\al,\al-a]$.

\noindent {\bf Definition of $h$ on $B(\hat w_{m,j};\epsilon)\cap[-\al,\al-a]$:}
On this interval, put $h=0$.
If
$1\leq m\leq \kappa$, then
$h$ is continuous at $y=\hat w_{m,j}$  and
\begin{equation}\label{b-10}
 \lim_{y\rightarrow w_{m,j}}
\left\{ h(y     -(1-ab)\frac{m}{b} )L_{m}(y) \right\}=0.
\end{equation}
 Hence the condition in Lemma \ref{b-3} (2) holds, $i.e.$, we have
now completed the proof of (a).


\noindent {\bf Definition of $h$ on $B(\hat w_{m,j}+a;\epsilon)\cap[a-\al,\al]$:}
We define $h$ on this set by
$ h(x)= \frac{b-g(x-a)h(x-a)}{g(x)}= \frac{b}{g(x)};$ thus $h$ is bounded here
by the choice of $\epsilon$ and \eqref{epsilon2-1}, and (b) holds
on $ B(\hat w_{m,j};\epsilon)\cap[-\al,\al-a].$

To summarize all these, let
$B_+:=\cup_{n=0}^{\kappa} \cup_{i=1}^{r_n}
\left(B(\tilde y_{n,i};\epsilon)\cap[a-\al,\al]\right),$
and $B_-:=\cup_{m=0}^{\kappa} \cup_{j=1}^{l_m}
\left(B(\hat w_{m,j};\epsilon)\cap [-\alpha,\al-a]\right).$
We have defined $h$ as a bounded function on $B:=[\al-a,a-\al] \cup
B_+ \cup (B_+ -a) \cup B_- \cup (B_- +a)$,
and (b) holds on
\bee \label{4101a} [-a,-\al]\cup[\al-a,0]\cup (B_+ -a)\cup B_-=[-a,-\al] \cup \left( B \cap [-a,0]\right).\ene


 \noindent {\bf Definition of $h$ on
$[-\al,\al]\setminus B$:}
Put $h=0$ on $\left([-\al,\al]\setminus B  \right)\cap [-a,0]$.
Note that the zeroset of $g$ within
$[-\al,\al]$ consists of  $Y_0$ and $W_0,$ so $g(x)\neq 0$ for
$x\in \overline{[-\al,\al]\setminus B};$ using the
continuity of g implies that $\inf_{x\in [-\al,\al]\setminus B} |g(x)|>0$.
We define $h(x)=\frac{b}{g(x)},\ x\in \left([-\al,\al]\setminus B  \right)\cap [0,a];$  thus,
we have now defined $h$ everywhere as a bounded function, and we
just need to complete the proof of (b). Since we have proved
(b) on the set in \eqref{4101a}, we just need to verify
(b) on the set $]- \al, 0]\setminus \left( B\cap [-a,0] \right).$ Note that $h$ vanishes on this set and that
\bes ]- \al, 0]\setminus \left( B\cap [-a,0] \right)
& = & \left([-\al,\al]\setminus B  \right)\cap [-a,0] \\
& = & ]-\al,\al-a[\setminus \left((B_+-a)\cup B_-\right) \\
& = & \left( \left([-\al,\al]\setminus B  \right)\cap [0,a]
\right) -a, \ens
where we used that $-\al\in B_-,\ \al\in B_+.$
 Thus, by the definition of $h$ on
$\left([-\al,\al]\setminus B  \right)\cap [0,a]$ (b) holds on $ ]-
\al, 0]\setminus \left( B\cap [-a,0] \right),$ as desired. \ep


\noindent{\bf Acknowledgments:}
This research was supported by Basic Science Research Program
through the National Research Foundation of Korea(NRF) funded by
the Ministry of Education(2013R1A1A2A10011922).
The authors thank the reviewer for suggesting that we provide more insight into the intuition behind Theorem \ref{2-2-1}; this motivated us to include
Example \ref{4108e}.

\end{document}